\documentclass[a4paper]{article}
\usepackage{amsmath,amssymb}
\usepackage[latin1]{inputenc}
\usepackage[dvips]{epsfig}
\usepackage[english,frenchb]{babel}

\newcommand{\UU}{\mathsf{U}}

\newcommand{\bord}{\mathsf{L}}
\newcommand{\Int}{\mathcal{I}}
\newcommand{\Y}{\mathsf{Y}}
\newcommand{\MPhi}{\boldsymbol{\Phi}}
\newcommand{\MTheta}{\boldsymbol{\Theta}}
\newcommand{\MPsi}{\boldsymbol{\Psi}}
\newcommand{\MDelta}{\boldsymbol{\Delta}}

\newcommand{\forme}{\boldsymbol{A}}
\newcommand{\PL}{\operatorname{PL}}
\newcommand{\sous}{\curvearrowleft}

\newcommand{\gch}{\epsfig{file=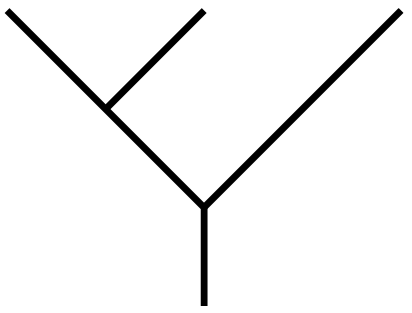,height=2mm}}
\newcommand{\drt}{\epsfig{file=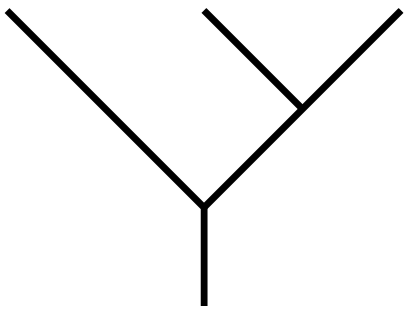,height=2mm}}
\newtheorem{theorem}{Théorème}[section] 
\newtheorem{proposition}[theorem]{Proposition}

\newtheorem{lemma}[theorem]{Lemme}

\newenvironment{proof}{\begin{trivlist}\item{\bf{Preuve.}}}
  {\hfill\rule{2mm}{2mm}\end{trivlist}}
\title{Sur le nombre d'intervalles\\ dans les treillis de Tamari}
\author{F. Chapoton}
\date{\today}

\begin{document}

\maketitle

\begin{abstract}
  On compte le nombre d'intervalles dans les treillis de Tamari. On
  utilise pour cela une description récursive de l'ensemble des
  intervalles. On introduit ensuite une notion d'intervalle nouveau
  dans les treillis de Tamari et on compte les intervalles nouveaux.
  On obtient aussi l'inverse de deux séries particulières dans un
  groupe de séries formelles en arbres.
\end{abstract}

\selectlanguage{english}

\begin{abstract}
  We enumerate the intervals in the Tamari lattices. For this, we
  introduce an inductive description of the intervals. Then a notion
  of ``new interval'' is defined and these are also enumerated. A a
  side result, the inverse of two special series is computed in a
  group of tree-indexed series.
\end{abstract}

\selectlanguage{frenchb}

\section{Introduction}

Les treillis de Tamari sont des ordres partiels remarquables, liés aux
polytopes de Stasheff (associaèdres). On peut notamment les réaliser
comme l'ordre induit par une forme linéaire particulière sur
l'ensemble des sommets d'une certaine version du polytope de Stasheff.
Ils ont été récemment généralisés par Reading \cite{reading}, qui a
introduit des treillis dits cambriens, associés aux carquois sur les
diagrammes de Dynkin. Dans ce contexte, les treillis de Tamari
correspondent aux carquois de type $\mathbb{A}$ équi-orientés.

L'objectif principal de cet article est de compter les intervalles
dans les treillis de Tamari. On y parvient en obtenant une description
récursive complète de l'ensemble des intervalles. On montre ainsi
l'existence d'une formule close pour le nombre d'intervalles.

Un objectif secondaire est de compter ceux parmi les intervalles qui
ne proviennent pas de treillis de Tamari d'indice inférieur par une
sorte d'``induction''. En termes géométriques, ce sont ceux qui ne
sont pas contenus dans une des facettes du polytope de Stasheff. On
les appelle les intervalles nouveaux, et on démontre aussi une formule
close pour le nombre d'intervalles nouveaux.

Il est assez remarquable que ces nombres d'intervalles et
d'intervalles nouveaux apparaissent aussi ensemble dans l'article
\cite{schaeffer} dont le sujet est l'énumération de certaines classes
de cartes planes.

Par ailleurs, on calcule l'inverse de deux éléments particuliers dans
un groupe de séries formelles en arbres lié à l'opérade dendriforme.
Ces séries sont des raffinements de la série génératrice usuelle des
nombres d'intervalles. La découverte de ces inverses a été la clé de
la description récursive des intervalles.

L'utilisation du système de calcul formel \textsf{MuPAD} a été
cruciale dans la recherche qui a abouti à cet article.

\section{Les treillis de Tamari}

Un \textbf{arbre binaire plan} est un graphe fini plan connexe et
simplement connexe, dont les sommets ont pour valence $1$ ou $3$, muni
d'un sommet de valence $1$ distingué appelé la \textbf{racine}. Les
autres sommets de valence $1$ sont appelés les \textbf{feuilles}, les
sommets de valence $3$ sont appelés \textbf{sommets internes}. On
dessine les arbres binaires plans avec les feuilles en haut et la
racine en bas.

Soit $Y_n$ l'ensemble des arbres binaires plans à $n$ sommets
internes. Le cardinal de $Y_n$ est le nombre de Catalan
$c_n=\frac{1}{n+1}\binom{2n}{n}$. Pour $n=1$, il y a un seul arbre
binaire plan, qui sera noté $\Y$. Sauf mention explicite du contraire,
on considère toujours qu'un arbre binaire plan a au moins un sommet
interne. L'utilisation éventuelle de l'arbre trivial à une feuille,
noté $|$, sera toujours explicitée.

\begin{figure}
  \begin{center}
    \scalebox{0.4}{\includegraphics{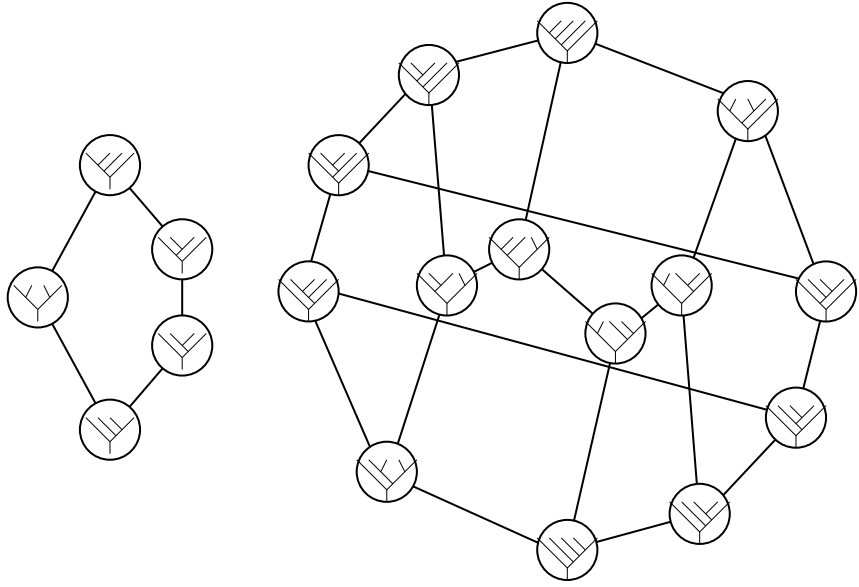}}
    \caption{Les treillis de Tamari $Y_3$ et $Y_4$}
    \label{tamari}
  \end{center}
\end{figure}

Le treillis de Tamari est un ordre partiel sur l'ensemble $Y_n$,
introduit par Tamari \cite{tamari}. Cet ordre partiel est défini comme
suit : un arbre $T$ est plus grand qu'un arbre $S$ si on peut passer
de $S$ à $T$ par une suite d'opérations consistant à remplacer
localement la configuration $\drt$ par la configuration $\gch$. Le
treillis $Y_1$ a un seul élément : $\Y$. Le treillis $Y_2$ est juste
$\drt \leq \gch$. La figure \ref{tamari} représente les diagrammes de
Hasse des treillis $Y_3$ et $Y_4$.

\subsection{Structures algébriques}

On définit des opérations $/$ et $\backslash$ sur les arbres binaires
plans (y compris l'arbre trivial $|$) : $S/T$ est obtenu en
identifiant la racine de $S$ avec la feuille la plus à gauche de $T$ ;
$S\backslash T$ est obtenu en identifiant la racine de $T$ avec la
feuille la plus à droite de $S$.  Ces opérations sont clairement
associatives. Par exemple $\Y/\Y=\gch$ et $\Y\backslash \Y=\drt$.

On considère l'espace vectoriel gradué $\mathcal{Y}$ ayant pour base
$Y_n$ en degré $n\geq 1$ et l'arbre trivial $|$ en degré $0$. On peut
étendre les opérations $/$ et $\backslash$ par linéarité.

Sur $\mathcal{Y}$, on dispose donc de deux produits associatifs. Un
troisième produit associatif existe, qui sera noté $*$, et qui peut
être défini comme suit :
\begin{equation}
  S * T =\sum_{S\backslash T \leq U \leq S/T} U,
\end{equation}
où la relation d'ordre est celle du treillis de Tamari. L'arbre
trivial $|$ est une unité pour les trois produits.

On va aussi utiliser l'opération $S \vee T$ définie par la greffe de
$S$ à gauche et de $T$ à droite sur les deux feuilles de $\Y$.

Pour plus de détails sur ces structures algébriques, on renvoie le
lecteur à \cite{arithm}.

\section{Énumération des intervalles}

Soit $\Int_n$ l'ensemble des intervalles dans le treillis de Tamari
$Y_n$.

En comptant les intervalles dans les premiers treillis de Tamari, on
obtient que les premiers termes de la suite $|\Int_n|$ sont, pour
$n\geq 1$,
\begin{equation}
  1,3,13,68,399,2530,etc.
\end{equation}

Une consultation de l'encyclopédie des suites \cite{oeis} mène
immédiatement à conjecturer le résultat suivant.
\begin{theorem}
  \label{theo_int}
  Le nombre d'intervalles dans le treillis de Tamari $Y_n$ est 
  \begin{equation}
    |\Int_n|=\frac{2(4n+1)!}{(n+1)!(3n+2)!}.
  \end{equation}
\end{theorem}

\begin{proof}
Il s'agit de calculer la série génératrice
\begin{equation}
  \phi=\sum_{n\geq 1} |\Int_n| y^n=y+3 y^2 + 13 y^3 + \dots
\end{equation}

Pour cela, on introduit un paramètre supplémentaire. Si $T$ est un
arbre binaire plan dans $Y_n$, on note $\bord(T)$ le nombre de
segments le long du bord gauche de $T$. Par exemple $\bord(\Y)=2$,
$\bord(\gch)=3$ et $\bord(\drt)=2$.

On considère alors la série génératrice raffinée
\begin{equation}
  \Phi=\sum_{n\geq 1}\sum_{[S,T]\in \Int_n} x^{\bord(T)} y^n=x^2 y+(2
  x^3+x^2) y^2 +(5 x^4 +5 x^3 + 3 x^2) y^3 + \dots,
\end{equation}
où la variable $x$ tient compte compte du paramètre $\bord$ pour le
maximum $T$ de l'intervalle $[S,T]$. Bien sûr, on retrouve $\phi$ en
posant $x=1$ dans $\Phi$.

On montre de façon combinatoire dans les sections \ref{decomposition},
\ref{description} et \ref{coeur} que $\Phi$ vérifie
\begin{equation}
  \Phi=x^2 y\left(1+\Phi/x\right)
\left(1+\left( \frac{\Phi-\phi}{x-1}\right)\right).
\end{equation}

En isolant $\phi$ dans un membre, on obtient
\begin{equation}
  \phi=\Phi+x-1+\frac{1}{y}\left( \frac{1}{x}-\frac{\Phi+1}{\Phi+x}\right),
\end{equation}
et comme $\phi$ ne dépend pas de $x$, le membre de droite est constant
en la variable $x$. Il en résulte que $\Phi$ est déterminée par
l'équation différentielle ordinaire
\begin{equation}
  \label{diffPhi}
  \partial_x \Phi=\frac{(x+\Phi)^2 (1/x^2-y) -1 -\Phi}{1-x+y(x+\Phi)^2}
\end{equation}
et la condition initiale $\Phi=0$ en $x=0$, car tout arbre binaire
plan a au moins deux segments sur son bord gauche.

On en déduit que $\Phi$ vérifie l'équation algébrique suivante :
\begin{multline}
\label{maxi8}
0={x}^{4}{y}^{4}\Phi^{8}+ ( 4\,{y}^{3}{x}^{3}+8\,{x}^{5}{y}^
{4}-{y}^{3}{x}^{4} ) \Phi^{7}+\\ ( -3\,{y}^{2}{x}^{3}+
32\,{y}^{3}{x}^{4}+6\,{y}^{2}{x}^{2}+28\,{x}^{6}{y}^{4}-3\,{x}^{5}{y}^
{3} ) \Phi^{6}+\\ ( 56\,{x}^{7}{y}^{4}+3\,{x}^{6}{y}^{
3}+ ( 108\,{y}^{3}-3\,{y}^{2} ) {x}^{5}-39\,{x}^{4}{y}^{2}+
40\,{y}^{2}{x}^{3}-3\,{x}^{2}y+4\,x y ) \Phi^{5}+\\ ( 1
+70\,{x}^{8}{y}^{4}+25\,{x}^{7}{y}^{3}+ ( 200\,{y}^{3}-9\,{y}^{2}
 ) {x}^{6}-146\,{x}^{5}{y}^{2}+\\ ( 116\,{y}^{2}+21\,y
 ) {x}^{4}-33\,{x}^{3}y+16\,{x}^{2}y-x ) \Phi^{4}+\\
 ( 56\,{x}^{9}{y}^{4}+45\,{x}^{8}{y}^{3}+ ( -6\,{y}^{2}+220
\,{y}^{3} ) {x}^{7}+ ( -254\,{y}^{2}-3\,y ) {x}^{6}+\\
 ( 184\,{y}^{2}+66\,y ) {x}^{5}-77\,{x}^{4}y+ ( -3+20
\,y ) {x}^{3}+3\,{x}^{2} ) \Phi^{3}+\\ ( 28\,{x}
^{10}{y}^{4}+39\,{x}^{9}{y}^{3}+ ( 6\,{y}^{2}+144\,{y}^{3}
 ) {x}^{8}+ ( -231\,{y}^{2}-5\,y ) {x}^{7}+\\ ( 
166\,{y}^{2}+68\,y ) {x}^{6}+ ( -3-67\,y ) {x}^{5}+
 ( 3+8\,y ) {x}^{4} ) \Phi^{2}+\\ ( 8\,{x}^
{11}{y}^{4}+17\,{x}^{10}{y}^{3}+ ( 52\,{y}^{3}+9\,{y}^{2}
 ) {x}^{9}+ ( -y-107\,{y}^{2} ) {x}^{8}+ ( -1+80
\,{y}^{2}+22\,y ) {x}^{7}+\\ ( 1-20\,y ) {x}^{6}
 ) \Phi+{x}^{12}{y}^{4}+3\,{x}^{11}{y}^{3}+ ( 3\,{y}^{
2}+8\,{y}^{3} ) {x}^{10}+ ( y-20\,{y}^{2} ) {x}^{9}+
 ( 16\,{y}^{2}-y ) {x}^{8}.
\end{multline}
En effet, la solution de cette équation algébrique vérifie l'équation
différentielle (\ref{diffPhi}) et la condition initiale.

On en déduit par spécialisation en $x=1$ (et simplification par
factorisation) que $\phi$ vérifie l'équation algébrique suivante :
\begin{multline}
  \label{eqphi}
  {y}^{3}{\phi}^{4}+ \left(4\,{y}+3
    \right)\,{y}^{2} {\phi}^{3}+ \left( 6\,{y}^{2}+17\,{y}+3 \right) y
    {\phi}^{2}\\+ \left( 4\,{y}^{3}+25\,{y} ^{2}-14\,y+1 \right)
    \phi+y\left({y}^{2}+11\,{y}-1\right)=0.
\end{multline}

Pour en déduire les coefficients de $\phi$, on peut utiliser comme
suit les résultats de \cite{schaeffer} dans le cas $a=b=1$. En posant
$\phi=\tau(1-\tau-\tau^2)$ dans (\ref{eqphi}), on obtient que $\tau$
vérifie $\tau=y(1+\tau)^4$. On calcule alors les coefficients de
$\phi$ par inversion de Lagrange. Ceci termine la démonstration du
Théorème.
\end{proof}

\section{Décomposition des intervalles}

\label{decomposition}

L'objet de cette section est l'obtention d'une description par
récurrence sur $n$ des intervalles.

\begin{lemma}
  Il existe une unique décomposition maximale de chaque arbre binaire
  plan $T$ en $T_1 / T_2 / \dots / T_k$.
\end{lemma}

\begin{proof}
  Clair, par découpe de chaque segment le long du bord gauche.
\end{proof}

Un arbre binaire plan $T$ est dit \textbf{indécomposable} si sa
décomposition maximale a un seul terme.

Étant donné un arbre binaire plan $T$, on définit une composition
$c(T)=(c_1,\dots,c_k)$ où $c_i$ est le nombre de sommets internes de
$T_i$ dans la décomposition maximale de $T$.

Soient $S$ et $T$ deux arbres binaires plans, et soient $(S_i)_i$ et
$(T_i)_i$ leurs décompositions maximales respectives.

\begin{lemma}
  \label{fusion}
  Si $S\leq T$, alors ou bien la composition $c(S)$ est strictement
  plus grossière que $c(T)$ ou bien $c(S)=c(T)$ et on a $S_i \leq
  T_i$ pour tout $i$.
\end{lemma}

\begin{proof}
  Il suffit de le montrer pour les relations élémentaires locales
  ($\drt \leq \gch$) qui engendrent l'ordre dans les treillis de
  Tamari.  Si cette relation est utilisée dans un des termes de la
  décomposition, la composition $c$ ne change pas et l'ordre est
  respecté dans chaque terme, avec changement seulement dans un terme.
  Si la relation est utilisée à la liaison entre deux termes, le long
  du bord gauche, alors il y a fusion des deux termes consécutifs
  impliqués en un seul terme, donc la composition devient strictement
  plus grossière et ne peut que le rester.
\end{proof}

\begin{figure}
  \begin{center}
    \scalebox{0.4}{\includegraphics{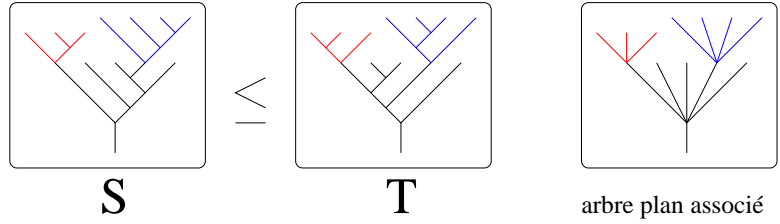}}
    \caption{Un intervalle dans $Y_9$}
    \label{exemple}
  \end{center}
\end{figure}

Par exemple, dans la figure \ref{exemple}, $c(S)=(2,7)$ et $c(T)=(1,1,2,4,1)$.

\begin{lemma}
  \label{factoriser_ordre}
  Soient $S=S_1/\dots/S_k$ et $T=T_1/\dots/T_k$ deux arbres binaires
  plans. On ne suppose pas que ces décompositions sont maximales. On
  suppose que $|S_i|=|T_i|$ pour tout $i$. Alors $S \leq T$ si et
  seulement si $S_i \leq T_i$ pour tout $i$.
\end{lemma}

\begin{proof}
  Clairement $S_i \leq T_i$ pour tout $i$ entraîne $S \leq T$.
  Montrons la réciproque. Choisissons une chaîne de relations d'ordre
  élémentaires de $S$ à $T$. Si une étape de cette chaîne se trouve le
  long du bord gauche à la liaison entre deux termes de la
  décomposition de $T$ fixée, la composition associée deviendrait
  strictement plus grossière par le lemme précédent, donc la
  composition associée à $S$ serait strictement plus grossière que
  celle de $T$, ce qui contredit l'hypothèse. Donc toutes les étapes
  de la chaîne se produisent à l'intérieur d'un des termes de la
  décomposition de $T$ fixée et la relation d'ordre se factorise.
\end{proof}

Étant donné deux intervalles $J$ et $K$, on peut définir un nouvel
intervalle $J/K$ avec $\min(J/K)=\min(J)/\min(K)$ et
$\max(J/K)=\max(J)/\max(K)$.

\begin{proposition}
  \label{decompo}
  Il existe une unique décomposition maximale de chaque intervalle $I$
  en $I_1 / I_2 / \dots / I_k$.
\end{proposition}

\begin{proof}
  Soit $S_1/\dots/S_k$ la décomposition maximale du minimum $S$ de
  $I$. Par le Lemme \ref{fusion}, la décomposition maximale du maximum
  $T$ de $I$ est plus fine que celle de $S$. Donc on peut en
  particulier décomposer le maximum $T$ comme le minimum $S$. Alors,
  chaque terme de la décomposition de $S$ est plus petit que chaque
  terme de $T$, par le Lemme \ref{factoriser_ordre}.

  Toute autre décomposition est un grossissement de celle-ci. En
  effet, elle est nécessairement un grossissement pour le minimum $S$.
\end{proof}

Dans l'exemple de la figure \ref{exemple}, la décomposition maximale a
deux termes, ayant respectivement $2$ et $7$ sommets internes.

On appelle \textbf{intervalle indécomposable} un intervalle dont la
décomposition maximale a un seul terme. On remarque que les
intervalles indécomposables sont exactement ceux dont le minimum est
indécomposable.

\section{Description des intervalles indécomposables}

\label{description}

Étant donné un intervalle $J$ et un segment $s$ du bord gauche de son
maximum $\max(J)$, on peut définir un nouvel intervalle $\Y *_s J$ où
$\max(\Y *_s J)$ est obtenu en rajoutant à gauche une arête sur
l'arête $s$ de $\max(J)$ et $\min(\Y *_s J)$ est obtenu en rajoutant à
gauche une arête sur l'arête racine de $\min(J)$. On peut aussi écrire
$\min(\Y *_s J)=\Y \backslash \min(J)$. Pour justifier que ceci est
bien un intervalle, on a les relations
\begin{equation}
  \Y \backslash \min(J) \leq \Y \backslash \max(J) \leq \max(\Y *_s J),
\end{equation}
où la seconde relation est obtenue par une suite de relations
élémentaires faisant remonter l'arête gauche le long du bord gauche de
$\max(J)$, jusqu'à atteindre le segment $s$.

\begin{figure}
  \begin{center}
    \scalebox{0.4}{\includegraphics{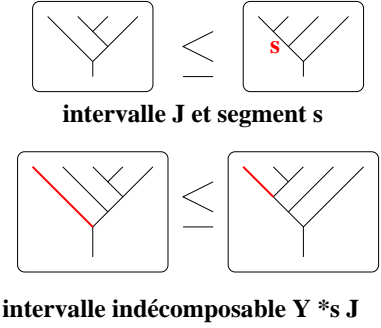}}
    \caption{Un intervalle indécomposable}
    \label{int_indec}
  \end{center}
\end{figure}

En particulier, l'élément $\min(\Y *_s
J)$ est indécomposable, donc l'intervalle $\Y *_s J$ est
indécomposable. Voir la figure \ref{int_indec} pour une illustration
de cette construction.

\begin{proposition}
  \label{indecompo}
  Tout intervalle indécomposable $I$ est ou bien $[\Y,\Y]$ ou bien de
  la forme $\Y *_s J$ pour un certain intervalle $J$ et un segment $s$
  du bord gauche de $\max(J)$, uniquement déterminés.
\end{proposition}

\begin{proof}
  Excluons le cas de l'intervalle $[\Y,\Y]$.

  Nécessairement le minimum de $I$ est de la forme $Y \backslash S$
  pour un certain arbre $S$. On peut aussi définir un arbre $T$ en
  enlevant l'arête gauche dans le maximum de $I$. 

  Montrons que $S \leq T$. Prenons une chaîne de relations
  élémentaires de $\min(I)$ à $\max(I)$. Dans cette chaîne, certaines
  étapes font intervenir l'arête la plus à gauche. Si on enlève ces
  étapes, on obtient une chaîne de relations élémentaires de $S$ à
  $T$.

  Donc on peut enlever l'arête la plus à gauche du minimum et du
  maximum de $I$ pour obtenir un nouvel intervalle $J=[S,T]$.

  Pour reconstituer de manière unique $I$, il suffit de se donner $J$
  et le segment $s$ sur le bord gauche de $\max(J)$ sur lequel se
  trouvait l'arête enlevée.
\end{proof}

\section{Équation fonctionnelle}

\label{coeur}

On déduit de la section précédente une équation fonctionnelle pour $\Phi$.

Soit $\theta$ la série génératrice des nombres d'intervalles
indécomposables :
\begin{equation}
  \theta=\sum_{n\geq 1}\sum_{{[S,T]\in \Int_n}\atop {\text{ indéc.}}}  y^n=y+2 y^2+8 y^3+41y^4+\dots,
\end{equation}
et soit $\Theta$ la série génératrice raffinée correspondante :
\begin{equation}
  \Theta=\sum_{n\geq 1}\sum_{{[S,T]\in \Int_n} \atop {\text{ indéc.}}} x^{\bord(T)} y^n=x^2 y+(x^3+x^2)y^2+(2 x^4+3 x^3+3 x^2) y^3 + \dots
\end{equation}

Alors on a les relations suivantes. Par la décomposition unique d'un
intervalle en intervalles indécomposables (Prop. \ref{decompo}), on a
\begin{equation}
  \Phi=\Theta+\Phi \, \Theta/x.
\end{equation}

Par la description des intervalles indécomposables (Prop.
\ref{indecompo}) et la définition de la fonction $\bord$, on a
\begin{equation}
  \Theta=x^2 y+ y \,\tilde{\Phi},
\end{equation}
où $\tilde{\Phi}$ est obtenu est remplaçant chaque $x^k$ dans $\Phi$ par
$x^2+\dots+x^{k+1}$. Ceci peut s'écrire
\begin{equation}
  \Theta=x^2 y+x^2 y \left( \frac{\Phi-\phi}{x-1}\right).
\end{equation}

On a donc la relation fonctionnelle suivante :
\begin{equation}
  \label{relaphi}
  \Phi=x^2 y\left(1+\Phi/x\right)
  \left(1+\left( \frac{\Phi-\phi}{x-1}\right)\right).
\end{equation}

\section{Séries formelles en arbres binaires plans}

Cette section fait intervenir explicitement des notions algébriques
plus sophistiquées : l'opérade dendriforme (due à Loday
\cite{lecturenotes}) et le groupe de séries formelles associé. Les
résultats de cette section ne sont pas utilisés dans le reste de
l'article. On renvoie le lecteur aux articles \cite{groupes,arithm}
pour la notion d'opérade et pour le groupe associé.

Soit $\MPhi$ la série génératrice des nombres d'intervalles selon leur
maximum :
\begin{equation}
  \MPhi=\sum_{n \geq 1} \sum_{ [S,T] \in \Int_n} T=\Y+2\,\gch+\drt+\dots
\end{equation}
et $\MTheta$ la série génératrice similaire des intervalles
indécomposables :
\begin{equation}
   \MTheta=\sum_{n \geq 1} \sum_{{[S, T]\in \Int_n}\atop
     {[S,T] \text{ indéc.}}} T=\Y+\gch+\drt+\dots  
\end{equation}

Par la décomposition unique en intervalles indécomposables (Prop.
\ref{decompo}), on a
\begin{equation}
  \MPhi=\MTheta+\MPhi/\MTheta.
\end{equation}

Par la description des intervalles indécomposables (Prop.
\ref{indecompo}), on a
\begin{equation}
  \qquad \MTheta=\Y+\Y*\MPhi.
\end{equation}
En effet, l'action du produit à gauche par $\Y$ est précisément de
faire la somme sur les arbres obtenus en rajoutant un arête à gauche
sur chaque segment du bord gauche.

On a obtenu la proposition ci-dessous.

\begin{proposition}
  On a l'équation fonctionnelle suivante :
  \begin{equation}
    \label{relaF}
    \MPhi=\Y+\MPhi/\Y+\Y*\MPhi+\MPhi/(\Y*\MPhi).
  \end{equation}
  Cette relation caractérise la série $\MPhi$.
\end{proposition}

Ceci a pour conséquence les propositions ci-dessous.

\begin{proposition}
  L'inverse de la série $\Y-\gch$ dans le groupe $\mathsf{G}_{Dend}$ des séries
  formelles en arbres binaires plans est la série $\Y+\gch \circ_1
  \MPhi$.
\end{proposition}

\begin{proof}
  Par la définition de la loi de groupe de $\mathsf{G}_{Dend}$, on doit montrer
  \begin{equation}
    (\Y) \circ_1(\Y+\gch\circ_1 \MPhi)-\big{(}\gch \circ_2(\Y+\gch \circ_1 \MPhi)\big{)} \circ_1(\Y+\gch \circ_1 \MPhi)=\Y.
  \end{equation}
  Par linéarité et comme $\Y$ est une unité pour les compositions
  $\circ$, ceci devient
  \begin{equation}
    \Y+\gch\circ_1 \MPhi-\gch-\gch\circ_1 (\gch \circ_1 \MPhi)-\gch\circ_2(\gch
    \circ_1 \MPhi)-\big{(}\gch \circ_2(\gch \circ_1 \MPhi)\big{)} \circ_1(\gch \circ_1 \MPhi)=\Y.
  \end{equation}
  Soit encore
  \begin{equation}
    \label{etape}
    \gch\circ_1 \MPhi=\gch+\gch\circ_1 (\gch \circ_1 \MPhi)+\gch\circ_2(\gch
    \circ_1 \MPhi)+\big{(}\gch \circ_2(\gch \circ_1 \MPhi)\big{)} \circ_1(\gch \circ_1 \MPhi).
  \end{equation}

  Comme on a
  \begin{equation}
    x \circ_1 (\gch \circ_1 y)=y /x,
  \end{equation}
  \begin{equation}
    x \circ_1 (y \circ_1 z)=(x \circ_1 y) \circ_1 z
  \end{equation}
  et
  \begin{equation}
    \gch\circ_2 (\gch \circ_1 x)=\gch\circ_1 (\Y * x),
  \end{equation}
  on peut réécrire l'équation (\ref{etape}) comme suit :
  \begin{equation}
    \gch\circ_1 \MPhi=\gch\circ_1 (\Y+\MPhi/\Y+\Y*\MPhi+\MPhi/(\Y*\MPhi)).
  \end{equation}
  Ceci résulte immédiatement de la relation fonctionnelle (\ref{relaF}).
\end{proof}

Considérons la série $\MDelta$ définie par
\begin{equation}
  \MDelta=\sum_{T_1,T_2} (-1)^{t_1+t_2} (t_1+1)\, T_1 \vee  T_2
  =\Y-(2\gch+\drt)+\dots,
\end{equation}
où les sommes portent cette fois sur les arbres binaires plans, y
compris l'arbre trivial $|$, et où on note $t_i$ le nombre de sommets
internes de $T_i$.

\begin{proposition}
  L'inverse de la série $\MPhi$ dans le groupe $\mathsf{G}_{Dend}$ des séries
  formelles en arbres binaires plans est la série $\MDelta$.
\end{proposition}

\begin{proof}
  Il s'agit de montrer que $\MPhi \MDelta=\Y$ dans le groupe
  $\mathsf{G}_{Dend}$.  Par la relation (\ref{relaF}), il suffit de
  calculer le produit
  \begin{equation}
    \left(\Y+\MPhi/\Y+\Y*\MPhi+\MPhi/(\Y*\MPhi)\right)  \MDelta.
  \end{equation}
  Comme le groupe $\mathsf{G}_{Dend}$ est contenu dans une algèbre
  associative et que la loi de groupe est linéaire à gauche, ceci vaut
  \begin{equation}
    \Y\MDelta+(\MPhi/\Y)\MDelta+(\Y*\MPhi)\MDelta+(\MPhi/(\Y*\MPhi))\MDelta.
  \end{equation}
  Posons $H=\MPhi\MDelta$. En utilisant la définition de la loi de
  groupe, on obtient la relation
  \begin{equation}
    \label{relaY}
    H=\MDelta+\sum_{k} (H_k /\Y)\circ_{k+1}\MDelta
    +(\Y*H)\circ_1 \MDelta+\sum_{k} (H_k/(\Y*H))\circ_{k+1}\MDelta,
  \end{equation}
  où on utilise une décomposition $H=\sum_{k} H_k$ en composantes homogènes.

  Cette relation caractérise uniquement la série $H$ par récurrence.
  Il suffit donc de montrer que $\Y$ vérifie (\ref{relaY}) pour
  pouvoir conclure que $H=\Y$. Il faut donc calculer 
  \begin{equation}
    \MDelta+(\gch)\circ_2 \MDelta+(\gch+\drt)\circ_1 \MDelta
    +(\Y/(\gch+\drt))\circ_2\MDelta.
  \end{equation}
  
  Ceci vaut
  \begin{multline}
  \sum_{T_1,T_2} (-1)^{t_1+t_2} (t_1+1)\, T_1 \vee  T_2+
  \sum_{T_1,T_2} (-1)^{t_1+t_2} (t_1+1)\, (\Y * T_1) \vee  T_2\\
  +\sum_{T_1,T_2} (-1)^{t_1+t_2} (t_1+1)\, (T_1 \vee  T_2) / \Y+
  \sum_{T_1,T_2} (-1)^{t_1+t_2} (t_1+1)\, T_1 \vee  (T_2* \Y)\\
  +\sum_{T_1,T_2} (-1)^{t_1+t_2} (t_1+1)\, ((\Y * T_1) \vee  T_2)/\Y+
  \sum_{T_1,T_2} (-1)^{t_1+t_2} (t_1+1)\, (\Y *T_1) \vee  (T_2* \Y).
  \end{multline}

  En utilisant le fait que
  \begin{equation}
    \sum_{T \not=|} T = \sum_{S} \Y * S = \sum_{S} S * \Y,
  \end{equation}
  on change les indices de sommation et on obtient
  \begin{multline}
  \sum_{T_1,T_2} (-1)^{t_1+t_2} (t_1+1)\, T_1 \vee  T_2
  -\sum_{T_1\not=|,T_2} (-1)^{t_1+t_2} (t_1)\, T_1 \vee  T_2\\
  +\sum_{T_1,T_2} (-1)^{t_1+t_2} (t_1+1)\, (T_1 \vee  T_2) / \Y
  -\sum_{T_1,T_2\not=|} (-1)^{t_1+t_2} (t_1+1)\, T_1 \vee  T_2\\
  -\sum_{T_1\not=|,T_2} (-1)^{t_1+t_2} (t_1)\, (T_1 \vee  T_2)/\Y+
  \sum_{T_1\not=|,T_2\not=|} (-1)^{t_1+t_2} (t_1)\, T_1 \vee T_2.
  \end{multline}
  On simplifie les termes par paires (1 \& 4), (2 \& 6) et (3 \& 5) :
   \begin{equation}
  \sum_{T_1} (-1)^{t_1} (t_1+1)\, T_1 /\Y
  -\sum_{T_1\not=|} (-1)^{t_1} (t_1)\, T_1 /\Y
  +\sum_{T_1,T_2} (-1)^{t_1+t_2} (T_1 \vee  T_2) / \Y.
  \end{equation}
  On simplifie à nouveau :
   \begin{equation}
     \Y+\sum_{T_1\not=|} (-1)^{t_1} (t_1)\, T_1 /\Y
  +\sum_{T_1\not=|} (-1)^{t_1+1} T_1 / \Y.
  \end{equation}
  Ceci donne $\Y$, ce qui termine la démonstration.
\end{proof}

\section{Découpage en intervalles nouveaux}

Un \textbf{arbre plan} est un graphe fini plan connexe et simplement
connexe, dont les sommets ont une valence différente de $2$, muni d'un
sommet de valence $1$ distingué appelé la \textbf{racine}. Les autres
sommets de valence $1$ sont appelés les \textbf{feuilles}, les sommets
de valence différente de $1$ sont appelés \textbf{sommets internes}.
On dessine les arbres plans avec les feuilles en haut et la racine en
bas. Soit $P_n$ l'ensemble des arbres plans à $n+1$ feuilles. Les
arbres binaires plans introduits précédemment sont en particulier des
arbres plans.

Étant donnés un arbre plan $T$ et, pour chaque sommet interne
$s$ de $T$ de valence $v(s)$, un intervalle $I_s$ dans $Y_{v(s)-2}$,
on peut définir un intervalle $G(T,(I_s)_s)$ en greffant les minimum
(resp.  les maximum) des intervalles $I_s$ selon le schéma de greffe
fourni par l'arbre plan $T$. La figure \ref{exemple} montre cette
construction pour un arbre plan à trois sommets et trois intervalles.

On dit qu'un intervalle est \textbf{nouveau} s'il ne peut pas s'écrire
ainsi de façon non triviale, c'est-à-dire s'il n'est pas de la forme
$G(T,(I_s)_s)$ avec $T$ ayant au moins deux sommets internes.

Remarque : tout intervalle nouveau est indécomposable. En effet, un
intervalle décomposable n'est évidemment à fortiori pas nouveau.

On appelle \textbf{découpage en nouveaux} de $I$ une écriture de $I$
sous la forme $G(T,(I_s)_s)$ où tous les $I_s$ sont nouveaux.
L'existence d'une telle écriture est claire par la définition des
intervalles nouveaux.

\begin{lemma}
  \label{descente}
  Soit $K$ un intervalle et $s$ un segment du bord gauche de $K$. Si
  $\Y *_s K$ est nouveau, alors tout découpage en nouveaux de $K$ est
  de la forme $K_1/ \dots / K_k$ où les $K_i$ sont nouveaux.
\end{lemma}

\begin{proof}
  Par l'absurde. Supposons qu'il existe un découpage de $K$ d'une
  autre forme. On peut donc en particulier décrire $K$ par la greffe
  de $K''$ sur une feuille (autre que la feuille gauche) de $K'$ pour
  certains intervalles $K'$ et $K''$. Mais on en déduit alors un
  découpage de $\Y *_s K$ sous la forme d'une greffe de $K''$ sur une
  feuille de $\Y *_s K'$. Ceci contredit l'hypothèse.
\end{proof}

\begin{proposition}
  \label{decoupage}
  Tout intervalle $I$ a un unique découpage en nouveaux.
\end{proposition}

\begin{proof}
  Il reste à montrer l'unicité. On procède par récurrence sur $n$, en
  utilisant la description des intervalles obtenue précédemment.

  Supposons d'abord que l'intervalle $I$ est décomposable. On choisit
  un découpage de $I$ en intervalles nouveaux. Ce découpage induit par
  regroupement une décomposition de $I$. Cette décomposition est
  nécessairement moins fine que la décomposition maximale.

  Si elle est strictement moins fine, un des termes admet une
  décomposition non triviale. Le découpage en nouveaux de $I$ induit
  un découpage en nouveaux de ce terme. L'intervalle situé à la racine
  de ce découpage en nouveaux est nécessairement aussi décomposable.
  Ce qui est absurde, car il est nouveau, donc indécomposable.
  
  Par conséquent, la décomposition induite par le découpage est égale
  à la décomposition maximale. Chacun des termes de la décomposition
  hérite alors d'un découpage en nouveaux, ces découpages sont uniques
  par récurrence, donc le découpage en nouveaux de $I$ est unique. Il
  s'obtient par recollement des découpages en nouveaux des termes de
  la décomposition maximale de $I$.

  \medskip

  \begin{figure}
    \begin{center}
      \scalebox{0.4}{\includegraphics{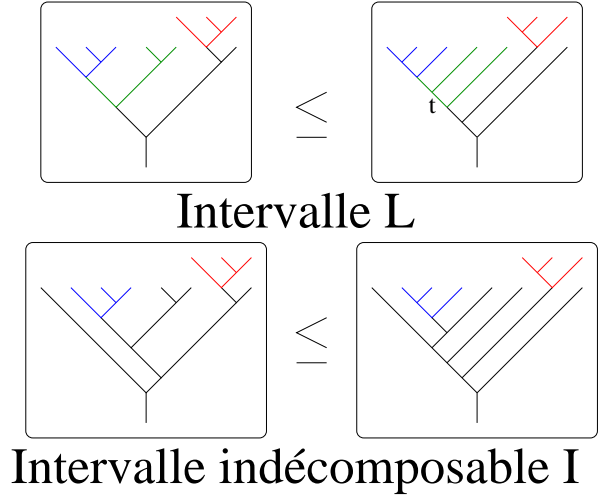}}
      \caption{Découpage en nouveaux d'un intervalle indécomposable}
      \label{fig-desc}
    \end{center}
  \end{figure}

  Supposons maintenant que l'intervalle $I$ est indécomposable. On
  peut supposer que $I$ n'est pas l'intervalle $[\Y,\Y]$. Il s'écrit
  donc $\Y *_t L$ pour certains $t$ et $L$. On choisit un découpage de
  $I$ en intervalles nouveaux. Soit $J$ l'intervalle nouveau situé à
  la racine de ce découpage.  Nécessairement, dans le découpage de
  $I$, il n'y a pas d'intervalle nouveau greffé sur la feuille de
  gauche de $J$, sinon $I$ serait décomposable.

  Comme $J$ est nouveau, donc indécomposable, $J$ s'écrit $\Y*_s K$
  pour certains $s,K$, ou bien $J=[\Y,\Y]$.

  Supposons d'abord que $J$ n'est pas $[\Y,\Y]$. Ce cas assez complexe
  est illustré par la figure \ref{fig-desc}. Par le Lemme
  \ref{descente}, le découpage en nouveaux de $K$ est de la forme $K_1
  / \dots /K_k$. On en déduit par recollement de ce découpage de $K$
  avec la partie supérieure du découpage en nouveaux de $I$ (les
  termes autres que $J$) un découpage en nouveaux de $L$.  Ce
  découpage en nouveaux de $L$ est unique par récurrence.

  Donc le découpage en nouveaux de $I$ est aussi unique et s'obtient à
  partir de celui de $L$ par recollement des morceaux $K_1,\dots,K_k$,
  qui sont les morceaux situés le long du bord gauche entre le segment
  racine et le segment $t$.

  Supposons maintenant que $J$ est l'intervalle $[\Y,\Y]$. Alors $I$
  est de la forme $\Y \backslash L$ pour un certain intervalle $L$, où
  l'opération $\backslash$ sur les intervalles est définie de manière
  similaire à l'opération $/$. Prenons un autre découpage en nouveaux
  de $I$ et soit $J'$ l'intervalle nouveau à la racine de ce
  découpage. Si $J'$ n'est pas l'intervalle trivial $[\Y,\Y]$, on peut
  facilement construire un découpage de $J'$, ce qui est absurde.

  Par conséquent, tout découpage en nouveaux de $I$ induit un
  découpage en nouveaux de $L$, unique par récurrence. Donc le
  découpage en nouveaux de $I$ est unique, obtenu par greffe du
  découpage en nouveaux de $L$ sur la feuille droite de l'intervalle
  $[\Y,\Y]$.
\end{proof}

L'exemple de la figure \ref{exemple} montre un intervalle et son
découpage en trois intervalles nouveaux, ainsi que l'arbre plan
associé à ce découpage.

Remarque : une autre démonstration de la Proposition \ref{decoupage}
peut sans doute être obtenue en utilisant la description géométrique
du treillis de Tamari par le biais d'une forme linéaire sur le
polytope de Stasheff.

\section{Résultat auxiliaire}

\label{muriel!}

On va utiliser ici un troisième type d'arbres. 

Un \textbf{arbre enraciné} est un graphe fini connexe et simplement
connexe muni d'un sommet distingué appelé la \textbf{racine}. On
dessine les arbres enracinés avec la racine en bas, en choisissant un
plongement arbitraire dans le plan. Par convention, on oriente les
arêtes de haut en bas.

Un automorphisme d'un arbre enraciné $A$ est une permutation de ses
sommets qui fixe la racine et préserve la relation d'adjacence. On
note $\sigma_A$ le cardinal du groupe d'automorphismes de $A$.

On considère l'espace vectoriel gradué $\PL$ ayant pour base en degré
$n$ les arbres enracinés à $n$ sommets.

Cet espace $\PL$ est muni d'un produit bilinéaire non-associatif : si
$A$ et $B$ sont des arbres enracinés, on pose
\begin{equation}
  A \sous B = \sum_{s \in A} A \sous_s B,
\end{equation}
où la somme porte sur les sommets de $A$ et $A \sous_s B$ est l'arbre
enraciné obtenu en ajoutant, dans l'union disjointe de $A$ et de $B$,
une arête entre la racine de $B$ et le sommet $s$ de $A$.

Soit $\nu$ une série formelle en une variable $y$. On associe à
chaque arbre enraciné $A$ une fonction $\nu_A$ définie
récursivement comme suit. On écrit $A$ comme la greffe d'arbres
enracinés $A_1,\dots,A_k$ (ensemble éventuellement vide) sur une
racine et on pose
\begin{equation}
  \nu_{A}= \left(\nu_{A_1} \dots \nu_{A_k}\right) \nu^{(k)},
\end{equation}
où $\nu^{(k)}$ désigne la $k$-ième dérivée de $\nu$ par rapport
à $y$. On a alors la propriété suivante.

\begin{proposition}
  Pour tous les arbres enracinés $A$ et $B$, on a
  \begin{equation}
    \nu_{A \sous B} = (\nu_A)' \nu_B.
  \end{equation}
\end{proposition}

\begin{proof}
  Ceci résulte de la théorie des algèbres pré-Lie, voir \cite{imrn}.
\end{proof}

Considérons maintenant la série
\begin{equation}
  \mathsf{\UU}=\sum_{A} \frac{1}{\sigma_A} A,
\end{equation}
qu'on peut voir comme appartenant à un complété de $\PL$. On décompose
$\UU$ selon le nombre de sommets dans les arbres :
\begin{equation}
  \mathsf{\UU}=\sum_{k \geq 0} \frac{\UU_{k+1}}{k!},
\end{equation}
où le terme $\UU_{k+1}$ est une combinaison linéaire d'arbres à $k+1$ sommets.

La proposition ci-dessous montre que la série formelle en arbres
enracinés $\UU$ s'identifie, aux signes près, à l'idempotent introduit
par Livernet \cite{muriel} dans l'étude des algèbres pré-Lie.

\bigskip

\begin{proposition}
  Les éléments $\UU_k$ sont déterminés par la récurrence
  \begin{equation}
  \UU_{k+1}
  =\sum_{\ell=1}^{k} \binom{k-1}{\ell-1} \UU_\ell \sous\UU_{k+1-\ell},    
  \end{equation}
  pour $k\geq 1$ et la condition initiale que $\UU_1$ est l'arbre
  enraciné à un sommet.
\end{proposition}

\begin{proof}
  Fixons un entier $k$ positif ou nul et un entier $\ell$ strictement
  positif et inférieur ou égal à $k$.

  Fixons un arbre enraciné $A$ à $k+1$ sommets, $B$ un arbre enraciné
  à $\ell$ sommets, $C$ un arbre enraciné à $k+1-\ell$ sommets,

  Considérons l'ensemble $E_1(A,B,C)$ des triplets $(s,i,j)$ où $s$
  est un sommet de $A$ distinct de la racine, $i$ un isomorphisme
  entre $C$ et le sous arbre de $A$ de racine $s$ et $j$ un
  isomorphisme entre $B$ et le complémentaire de ce sous-arbre de $A$.

  Alors, cet ensemble $E_1(A,B,C)$ est en bijection avec l'ensemble
  $E_2(A,B,C)$ des couples $(s,\varepsilon)$, où $s$ est un sommet de
  $B$ et $\varepsilon$ un isomorphisme entre $A$ et l'arbre $B \sous_s
  C$. La bijection utilise les isomorphismes pour identifier le sommet
  $s$ de $A$ à un sommet $s$ de $B$. On déduit l'isomorphisme
  $\varepsilon $ des deux isomorphismes $i$ et $j$ et réciproquement.

  On a donc une égalité $\sum_{E_1(A,B,C)} A = \sum_{E_2(A,B,C)} B \sous_s C $.

  De plus, les groupes d'automorphismes de $B$ et $C$ agissent
  librement sur $E_1(A,B,C)$ et le groupe d'automorphisme de $A$ agit
  librement sur $E_2(A,B,C)$.

  Calculons maintenant 
  \begin{align}
    \frac{\UU_{k+1}}{k!}&=\sum_{A} \frac{1}{\sigma_A} A
    = \frac{1}{k} \sum_{(A,s\in A)} \frac{1}{\sigma_A} A\\
    &=\frac{1}{k} \sum_{(A,B,C)} \sum_{(s,i,j) \in E_1(A,B,C)} \frac{1}{\sigma_A \sigma_B \sigma_C} A\\
    &=\frac{1}{k} \sum_{(A,B,C)} \sum_{(s,\varepsilon)\in E_2(A,B,C)} \frac{1}{\sigma_A \sigma_B \sigma_C} B \sous_s C\\
      &=\frac{1}{k} \sum_{(B,C,s \in B)} \frac{1}{\sigma_B \sigma_C} B \sous_s C\\
      &=\frac{1}{k} \sum_{(B,C)} \frac{1}{\sigma_B \sigma_C} B \sous C\\
      &=\frac{1}{k} \left(\sum_B \frac{1}{\sigma_B} B\right) \sous \left(\sum_C
      \frac{1}{\sigma_C} C\right)\\
      &= \frac{1}{k} \frac{\UU_{\ell}}{(\ell-1)!} \sous
      \frac{\UU_{k+1-\ell}}{(k-\ell)!}.
  \end{align}
  Ceci termine la démonstration de la proposition.
\end{proof}

\section{Énumération des intervalles nouveaux}

Si on compte les intervalles nouveaux dans les premiers treillis de
Tamari, on trouve
\begin{equation}
  1,1,3,12,56,288,1584,etc,
\end{equation}
pour $n \geq 1$.

On est mené, par consultation de l'encyclopédie des suites d'entiers,
à conjecturer le résultat suivant.

\begin{theorem}
  \label{theo_int_nouv}
  Pour $n \geq 2$, le nombre d'intervalles nouveaux dans le treillis de
  Tamari $Y_n$ est
  \begin{equation}
    3 \,\frac{ 2^{n-2} (2n-2)!}{(n-1)!(n+1)!}.
  \end{equation}
\end{theorem}

La preuve est l'objet de la suite de cette section.

\begin{figure}
  \begin{center}
    \scalebox{0.4}{\includegraphics{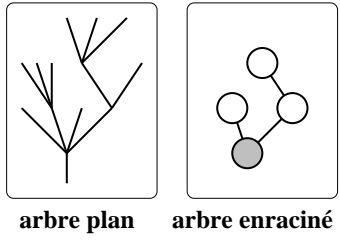}}
    \caption{Un arbre plan et l'arbre enraciné associé}
    \label{forme}
  \end{center}
\end{figure}

A chaque arbre plan $T$, on peut associer un arbre enraciné
$\forme(T)$ (la \textbf{forme} de $T$) de la façon suivante : les
sommets de $\forme(T)$ sont les sommets internes de $T$, les arêtes de
$\forme(T)$ sont les arêtes entre sommets internes de $T$, la racine
de $\forme(T)$ est le sommet interne de $T$ qui est adjacent à la
racine de $T$, voir la figure \ref{forme}.

\medskip

On introduit par commodité une version décalée de la série génératrice
des intervalles :
\begin{equation}
  \psi=y\phi=\sum_{n\geq 1} |\Int_n| y^{n+1}=y^2+3 y^3 + 13 y^4 + \dots
\end{equation}
et une série génératrice similaire pour les intervalles nouveaux :
\begin{equation}
  \nu=\sum_{n\geq 1} \sum_{{[S,T]\in\Int_n}\atop{\text{nouveau}}} y^{n+1}=y^2+
    y^3 + 3 y^4 + 12 y^5 +\dots
\end{equation}
On peut déduire de l'équation algébrique (\ref{eqphi}) pour $\phi$ une équation
algébrique pour $\psi$ : 
\begin{multline}
    {\psi}^{4}+ \left( 4\,y+3 \right) {\psi}
^{3}+ \left( 6\,{y}^{2}+17\,y +3\right) {\psi}^{2} \\+ \left( 
4\,{y}^{3}+25\,{y}^{2}-14\,y+1 \right) \psi+y^2
\left({y}^{2}+11\,{y}-1 \right) =0
\end{multline}

On va chercher à obtenir une équation algébrique pour $\nu$.

\medskip

Soit $N_n$ le nombre d'intervalles nouveaux dans le treillis de Tamari
$Y_n$. Le découpage unique d'un intervalle en intervalles nouveaux
donne la relation suivante :
\begin{equation}
  \label{nouveau_brut}
  \psi=\sum_{n \geq 1} \sum_{T\in P_n} N_T \,y^{n+1},
\end{equation}
où $P_n$ est l'ensemble des arbres plans à $n+1$ feuilles et $N_T$ est
le produit sur l'ensemble des sommets internes $s$ de $T$ des
$N_{v(s)-2}$ ($v(s)$ est la valence de $s$).

Dans cette somme, on regroupe les termes selon la valeur de
$\forme(T)$ :
\begin{equation}
  \psi=\sum_{A}\sum_{n\geq 1} \sum_{{T \in P_n}\atop{\forme(T)\simeq A}}
  N_T \, y^{n+1},
\end{equation}
où la somme externe porte sur les arbres enracinés.

Fixons un arbre enraciné $A$ et soit $I$ l'ensemble de ses sommets.
Notons $v_i$ le nombre d'arêtes entrantes du sommet $i$. La somme des
$v_i$ vaut $|I|-1$, car chaque sommet sauf la racine a exactement une
arête sortante.

On a besoin de décrire l'ensemble des arbres plans de forme isomorphe
à $A$. On a un bijection entre d'une part l'ensemble $F_1(A)$ formé
par les paires $(T,\gamma)$ où $T$ est un arbre plan et $\gamma$ un
isomorphisme entre $\forme(T)$ et $A$ et d'autre part l'ensemble
$F_2(A)$ formé par les paires $((\ell_i)_i,(Z_i)_i)$ où, pour chaque
sommet $i$ de $A$,
\begin{itemize}
\item $\ell_i$ est un nombre supérieur ou égal à $2$,
\item $Z_i$ est une injection de l'ensemble des arêtes entrantes en $i$
  dans l'ensemble $\{1,\dots,\ell_i\}$.
\end{itemize}

La bijection est la suivante : à la paire $(T,\gamma)$, on associe les
nombres $\ell_i$ tels que $\ell_i+1$ est la valence du sommet de $T$
correspondant par $\gamma$ au sommet $i$ de $A$. Par $\gamma$, les
arêtes entrantes en $i$ dans $A$ sont envoyées dans l'ensemble des
arêtes entrantes du sommet correspondant de $T$, ce qui donne
l'injection $Z_i$ voulue.

Réciproquement, si on connaît les nombres $\ell_i$ et les injections
$Z_i$, on peut aisément reconstruire un arbre plan $T$ muni d'une
bijection $\gamma$.

% Un arbre plan $T$ de forme $A$ est uniquement déterminé par les
% données suivantes :
% \begin{itemize}
% \item un dessin de $A$ dans le plan : un ordre total sur les arêtes
%   entrantes de chaque sommet de $A$, modulo l'action libre du groupe
%   des automorphismes de $A$.
% \item les valences $\ell_i+1$ dans $T$ de ses sommets internes (en
%   bijection avec les sommets de $A$)
% \item les positions de greffe : pour chaque $i$, $v_i$
%   entiers distincts dans $\{1,\dots,\ell_i\}$.
% \end{itemize}

La contribution de $A$ à $\psi$ est donc donnée par
\begin{equation}
  \frac{1}{\sigma_A} \sum_{(\ell_i)_{i \in I}} \prod_{i \in I} \left( \frac{\ell_i
      !}{(\ell_i-v_i)!} N_{\ell_i-1} \right) y^{n+1},
\end{equation}
où $n$ est la somme des $\ell_i-1$. En répartissant la variable $y$,
ceci vaut encore
\begin{equation}
  \frac{1}{\sigma_A} \sum_{(\ell_i)_{i \in I}} \prod_{i \in I} \left(\frac{\ell_i
    !}{(\ell_i-v_i)!} N_{\ell_i-1} y^{\ell_i-v_i}\right),
\end{equation}
soit
\begin{equation}
  \frac{1}{\sigma_A} \prod_{i\in I} \left(\sum_{\ell_i \geq 2} \frac{\ell_i
    !}{(\ell_i-v_i)!}N_{\ell_i-1} y^{\ell_i-v_i}\right).
\end{equation}
On reconnaît dans les facteurs de cette formule les dérivées de $\nu$,
on obtient donc
\begin{equation}
  \frac{1}{\sigma_A} \prod_{i\in I} \nu^{(v_i)}.
\end{equation}

En utilisant la définition des fonctions $\nu_A$ comme dans la section
\ref{muriel!}, on peut donc écrire
\begin{equation}
  \label{sommarb}
  \psi=\sum_{A} \frac{1}{\sigma_A} \nu_A=\nu+\nu \nu'+ \nu (\nu')^2 +\frac{1}{2} \nu^2 \nu''+ \dots
\end{equation}
où la somme porte sur les arbres enracinés, 

\medskip

Introduisons une suite de fonctions $(\alpha_k)_{k \geq 1}$ de la
variable $y$. Le terme initial est $\alpha_1=\nu$. Les termes suivants
sont donnés pour $k\geq 1$ par la récurrence
\begin{equation}
  \alpha_{k+1}
  =\sum_{\ell=1}^{k} \binom{k-1}{\ell-1}\partial_y(\alpha_\ell)\,\alpha_{k+1-\ell}.
\end{equation}

Alors, par les résultats de la section \ref{muriel!}, l'équation
(\ref{sommarb}) équivaut à
\begin{equation}
  \psi=\sum_{k \geq 0} \frac{\alpha_{k+1}}{k!}.
\end{equation}

Introduisons la série génératrice double
\begin{equation}
  \MPsi=\sum_{k \geq 0}{\alpha_{k+1}}\frac{z^k}{k!}.
\end{equation}

La récurrence définissant les fonctions $\alpha_k$ se traduit en une
équation aux dérivées partielles pour $\MPsi$ :
\begin{equation}
  \label{edp}
  \partial_z \MPsi= (\partial_y \MPsi) \MPsi,
\end{equation}
avec la condition initiale
\begin{equation}
  \label{condi}
  \MPsi\vert_{z=1}=\psi.
\end{equation}

Par conséquent, $\Psi$ est solution de l´équation algébrique
\begin{multline}
 {z}^{4}{\MPsi}^{4}+ \left( 4\,{z} y-8+11\,{z} \right) z^2 {\MPsi}^{3}+ \left( 6\,{z}^{2}{y}^{2}-{z}^
{2}+33\,{z}^{2}y-16\,z y+16-12\,z \right) {\MPsi}^{2}\\+ \left( -8\,{
y}^{2}-2\,z y+1+4\,z{y}^{3}-12\,y+33\,z{y}^{2} \right) \MPsi+11\,{y
}^{3}+{y}^{4}-{y}^{2} =0.
\end{multline}
En effet, la solution de cette équation algébrique vérifie l'équation
différentielle (\ref{edp}) et la condition initiale (\ref{condi}) en
$z=1$. En utilisant la spécialisation
\begin{equation}
  \MPsi\vert_{z=0}=\alpha_1=\nu,
\end{equation}
on déduit de cette équation algébrique pour $\MPsi$ une équation
algébrique pour $\nu$ :
\begin{equation}
  {y}^{2}-11\,{y}^{3}-{y}^{4}+ \left( -1+12\,y+8\,{y}^{2} \right) \nu
  -16\,  \nu^{2}=0.
\end{equation}
Donc
\begin{equation}
 \nu= \frac{1}{32}\left(-1+12\,y+8\,{y}^{2}+\,(1-8\,y)^{3/2}\right).
\end{equation}
On en déduit immédiatement la formule attendue pour les coefficients
de $\nu$. Ceci démontre le Théorème \ref{theo_int_nouv}.

\bibliographystyle{plain}
\bibliography{comptage}

\begin{thebibliography}{1}

\bibitem{imrn}
Fr{\'e}d{\'e}ric Chapoton and Muriel Livernet.
\newblock Pre-{L}ie algebras and the rooted trees operad.
\newblock {\em Internat. Math. Res. Notices}, (8):395--408, 2001.

\bibitem{groupes}
Fr{é}d{é}ric Chapoton.
\newblock {Rooted trees and an exponential-like series}.
\newblock arXiv:math.QA/0209104.

\bibitem{schaeffer}
Robert Cori and Gilles Schaeffer.
\newblock Description trees and {T}utte formulas.
\newblock {\em Theoret. Comput. Sci.}, 292(1):165--183, 2003.
\newblock Selected papers in honor of Jean Berstel.

\bibitem{tamari}
Samuel Huang and Dov Tamari.
\newblock Problems of associativity: {A} simple proof for the lattice property
  of systems ordered by a semi-associative law.
\newblock {\em J. Combinatorial Theory Ser. A}, 13:7--13, 1972.

\bibitem{muriel}
Muriel Livernet.
\newblock {A rigidity theorem for preLie algebras}.
\newblock {à} para{î}tre dans JPAA.

\bibitem{lecturenotes}
Jean-Louis Loday.
\newblock Dialgebras.
\newblock In {\em Dialgebras and related operads}, volume 1763 of {\em Lecture
  Notes in Math.}, pages 7--66. Springer, Berlin, 2001.

\bibitem{arithm}
Jean-Louis Loday.
\newblock Arithmetree.
\newblock {\em J. Algebra}, 258(1):275--309, 2002.
\newblock Special issue in celebration of Claudio Procesi's 60th birthday.

\bibitem{reading}
Nathan Reading.
\newblock {Cambrian Lattices}.
\newblock arXiv:math.CO/0402086.

\bibitem{oeis}
N.~J.~A. Sloane.
\newblock The {O}n-{L}ine {E}ncyclopedia of {I}nteger {S}equences, 2006.
\newblock www.research.att.com/{~}njas/sequences/.

\end{thebibliography}

\end{document}